\documentclass[a4paper,11pt]{article}
\usepackage[T1]{fontenc}
\usepackage[utf8]{inputenc}
\usepackage[english]{babel}
\usepackage{lmodern}
\usepackage{amsmath,amssymb,amsthm,mathtools}
\usepackage{geometry}
\usepackage{booktabs,array,tabularx,float}
\usepackage{tikz}
\usepackage[hidelinks]{hyperref}

\newtheorem{theorem}{Theorem}[section]
\newtheorem{lemma}[theorem]{Lemma}
\newtheorem{proposition}[theorem]{Proposition}
\newtheorem{corollary}[theorem]{Corollary}
\theoremstyle{remark}

\newcommand{\bbone}{\mathbf 1}
\newcommand{\Iarr}[2]{\{#1\,;\,#2\}}
\newcommand{\yes}{\mathrm{yes}}
\newcommand{\no}{\mathrm{no}}

\title{Boolean functions and completely regular codes of bounded degree on the halved cube}

\author{%
  Ivan Mogilnykh%
  \thanks{Sobolev Institute of Mathematics, Novosibirsk, Russia,
    \texttt{ivmog84@gmail.com} (corresponding author)}%
  \kern.35em\relax%
  \footnotemark[3]%
  \and
  Alexandr Valyuzhenich%
  \thanks{Sobolev Institute of Mathematics, Novosibirsk, Russia,
    \texttt{graphkiper@mail.ru}}%
  \kern.35em\relax%
  \thanks{The work was carried out within the state assignment of the
    Sobolev Institute of Mathematics SB RAS, project no. FWNF-2026-0011.}%
}

\date{}
\begin{document}
\maketitle

\begin{abstract}

We establish a bijection between Boolean functions of degree $d$ on the $n$-dimensional hypercube and Boolean functions of degree $d$ on the halved $n$-cube for $n\geq 4d+1$.
Moreover, for each fixed $d$ and sufficiently large $n$, we show that 
every completely regular code of degree $d$ in the halved $n$-cube
has covering radius $1$.

\end{abstract}  

\section{Introduction}
A function $f:\{0,1\}^n\rightarrow \{0,1\}$ is called {\em Boolean}.
Given a function $f$ on $\{0,1\}^n$, a variable $x_i$, $1\leq i\leq n$, is called {\em relevant}
if there are two vectors $x,y\in \{0,1\}^n$ such that $x$ and $y$ differ only in position $i$ and $f(x)\neq f(y)$.
The {\em $n$-dimensional hypercube} $H_n$ is defined as follows.
The vertex set of $H_n$ is $\{0,1\}^n$, and two vertices are adjacent if they differ in exactly one coordinate. The vertices of the {\em halved $n$-cube} $Q_n$ are all even-weight binary vectors of length $n$, 
and two vertices are adjacent if they differ in exactly two coordinates.

The famous Nisan--Szegedy theorem \cite{NS94} states that any degree-$d$ Boolean function on $H_n$ has at most $d\cdot 2^{d-1}$ relevant variables.
In 2020, Chiarelli, Hatami and Saks \cite{CHS20} improved this bound to $6.614\cdot 2^{d}$.
Subsequently, in 2022, Wellens \cite{W22} further improved it to $4.394\cdot 2^{d}$.
In Johnson graphs, Filmus and Ihringer \cite{FI19J} established an analogue
of the Nisan--Szegedy theorem for Boolean functions, showing that functions of fixed
degree are "juntas"\, when the parameters are sufficiently
large.

Now we briefly discuss the results on Boolean functions of degree at most $3$ on $H_n$.
Any degree-$1$ Boolean function on $H_n$ has exactly one relevant variable (this can be proved by induction on $n$; it also follows directly from the Nisan--Szegedy theorem).
By the Nisan--Szegedy theorem, any degree-$2$ Boolean function on $H_n$ has at most $4$ relevant variables.
Therefore, all such functions can be found by the enumeration of all Boolean functions on $H_4$.
In other terms, such functions were characterized by Camion et al. in \cite{CCCS92}.
The results of Tarannikov and Kirienko \cite{TK00} on correlation-immune functions imply that
any degree-$3$ Boolean function on $H_n$ has at most $10$ relevant variables (this bound is sharp).
Based on this bound, in 2023 Rasoolzadeh \cite{R23} classified all $(n-4)$-resilient Boolean functions with $n$ variables.
The classification of degree-$3$ Boolean functions can be derived from Rasoolzadeh's classification (see details in \cite[Section 2.1]{KV24}).

Characteristic functions of completely regular codes (CRCs) whose
largest eigenvalue index is $d$ form an important class of Boolean
functions of degree $d$. Completely regular codes are also of
independent interest \cite{SS}. In particular, a Boolean function on
a distance-regular graph has degree $1$ if and only if its support
is a CRC with covering radius $1$ whose non-principal eigenvalue
has index $1$. 
Boolean functions of degree $1$ have been studied in classical
association schemes \cite{FI19,Ihringer24,Filmus26}.

Completely regular codes of covering radius $1$ whose characteristic
functions have degree $2$ have been investigated in Hamming graphs
\cite{MV20,MTV25} and Johnson graphs \cite{EGGV25,GG13,V20}.
Completely regular codes of covering radius $1$ in halved cubes were
studied in \cite{KrotovMogilnykhVasileva2018}, where all possible
parameter sets for such codes in $Q_n$ with $n\leq 10$ were determined.

In this work, we investigate Boolean functions of bounded degree on the
halved $n$-cube $Q_n$ and their connection to completely regular codes.
In Section~3, we establish a rigidity property: if $n\geq 4d+1$,
restriction to $Q_n$ gives a bijection between Boolean functions of
degree $d$ on the hypercube $H_n$ and Boolean functions of degree $d$
on $Q_n$. Moreover, this bound is sharp.

Having established the restriction result for Boolean functions of fixed
degree, we turn to completely regular codes in Section~4. In general, a CRC of degree $d$ has covering radius at most $d$ because 
a CRC with covering radius $\rho$ has exactly $\rho+1$ distinct
eigenvalues. For halved
cubes, however, we show that for every fixed $d$, all nontrivial
CRCs of degree $d$ in $Q_n$ have covering radius $1$ whenever $n$ is
sufficiently large. Moreover, their characteristic functions are
restrictions of characteristic functions of CRCs with covering radius
$1$ in the hypercube $H_n$.

Thus, for every fixed degree $d$, all nontrivial CRCs in sufficiently halved cubes of sufficiently large diameter have covering radius $1$. In contrast, Hamming, Johnson, and Grassmann graphs admit CRCs of degree $d$ and covering radius $d$ for every diameter at least $d$ (see Remark~1 in Section~4).

Using the correspondence established in Section~3 between Boolean
functions of fixed degree on the hypercube $H_n$ and on the halved
$n$-cube $Q_n$, we obtain a complete classification of Boolean functions
of degree at most $2$ on $Q_n$ for $n\geq 9$ in Section~5.
We now proceed to the preliminaries on Boolean functions.
The definitions related to CRCs are given in Section~4.

\section{Preliminaries and notation}

Let $\Gamma$ be a distance-regular graph of diameter $r$.
We index its eigenvalues in decreasing order:
\[
\theta_0(\Gamma)>\theta_1(\Gamma)>\cdots>\theta_r(\Gamma).
\]
Let $U_i(\Gamma)$ be the eigenspace corresponding to $\theta_i(\Gamma)$,
and $E_i(\Gamma)$ the orthogonal projector onto it.
A real-valued function $f$ on $V(\Gamma)$ has degree  $d$ if
\[
 f\in\bigoplus_{i=0}^{d}U_i(\Gamma),
\]
and $E_d(\Gamma)f\ne0$. The zero function is defined to have degree $0$.

Write $Z_2^n$ for the binary words of length $n$, and $|x|$
for the Hamming weight of $x$. The $n$-dimensional \emph{hypercube} 
$H_n$ has vertex set $V(H_n)= Z_2^n$, with two vertices adjacent
if they differ in exactly one position. The vertex set of $Q_n$ is
\[
 A_n=\{x\in Z_2^n:|x|\equiv0\pmod2\}.
\]

For a function $f$ on $H_n$, write $f^+:=f|_{A_n}$. The eigenvalues of the hypercube and halved cube graphs $H_n$ and $Q_n$ are
\[
 \theta_i(H_n)=n-2i\quad(0\le i\le n),\qquad
 \theta_i(Q_n)= \frac{(n-2i)^2-n}{2}
 \quad\bigl(0\le i\le\lfloor n/2\rfloor\bigr).
\]
For $u,x\in Z_2^n$, define the Fourier character
\[
 \chi_u(x)=(-1)^{\langle u,x\rangle},
\]
where the inner product is taken over $Z_2$. Write
$\overline u=u+(1,\ldots,1)$ for the bitwise complement of $u$.
Since both graphs are Cayley graphs of abelian groups,
\[
 U_i(H_n)=\operatorname{span}\{\chi_u:|u|=i\},
\]
and, for $0\le i<n/2$,
\[
 U_i(Q_n)=\operatorname{span}\{\chi_u^+:|u|=i\},
\]
where the displayed restrictions form a basis. When $n$ is even
and $i=n/2$, one takes one representative from each pair
$\{u,{\overline u}\}$, since
$\chi_{\overline u}^+=\chi_u^+$.
Thus, for $0\le i<n/2$, restriction is an isomorphism from each of
$U_i(H_n)$ and $U_{n-i}(H_n)$ onto $U_i(Q_n)$:
\begin{equation}\label{e}
 U_i(H_n)^+=U_{n-i}(H_n)^+=U_i(Q_n).
\end{equation}

\section{Boolean functions of bounded degree on $Q_n$}

\begin{lemma}\label{lem:extension}
If $n\ge2d+1$ then for every function $g$ of degree $d$
on $Q_n$ there is a unique function $f$ of degree $d$
on $H_n$ such that $f^+=g$. 
\end{lemma}
\begin{proof}
By~\eqref{e}, restriction is a bijection between
$\bigoplus_{i=0}^d U_i(H_n)$ and $\bigoplus_{i=0}^d U_i(Q_n)$.
\end{proof}

\begin{theorem}\label{cor:boolean-lifting}
For $n\ge4d+1$, restriction $f\mapsto f^+$ is a bijection
between Boolean functions of degree $d$ on $H_n$ and $Q_n$.
\end{theorem}

\begin{proof}
By Lemma~\ref{lem:extension}, every Boolean function $g$
of degree $d$ on $Q_n$ has a unique real-valued
extension $f$ of degree $d$ on $H_n$.

Consider the function $f^2-f$, which has degree
$i$ for some $0\le i\le 2d$, since the product
of two Fourier characters of weights at most $d$
has weight at most $2d$. Since $n\ge 4d+1\ge 2i+1$, Lemma~\ref{lem:extension}
implies that $(f^2-f)^+$ has degree $i$.
However, $(f^2-f)^+\equiv 0$, since $g=f^+$ is Boolean.
Therefore $i=0$ and $f^2-f\equiv 0$,
so $f$ is Boolean.
\end{proof}

Note that the  above bound  is sharp.

\begin{proposition}\label{prop:sharp}
For every  $d\ge1$, there exists a 
function $f$ of degree $d$ on $H_{4d}$ such that $f^+$ is
Boolean on $Q_{4d}$ but $f$ is not Boolean on $H_{4d}$.
\end{proposition}
\begin{proof}
The function 
$ \frac{x_1+x_2+x_3-x_4}{2}$ is of degree $1$ and Boolean on $Q_4$. However, if
$x=(1,1,1,0)$, the value of the function is $3/2$, so it is not Boolean on $H_n$. 

We replace each $x_i$ with degree $d$ function $\sum\limits_{j=1}^{d} x_{d(i-1)+j} \mbox{ mod }2$ to obtain the function on $H_{4d}$ of degree $d$ fulfilling the required properties.
\end{proof}

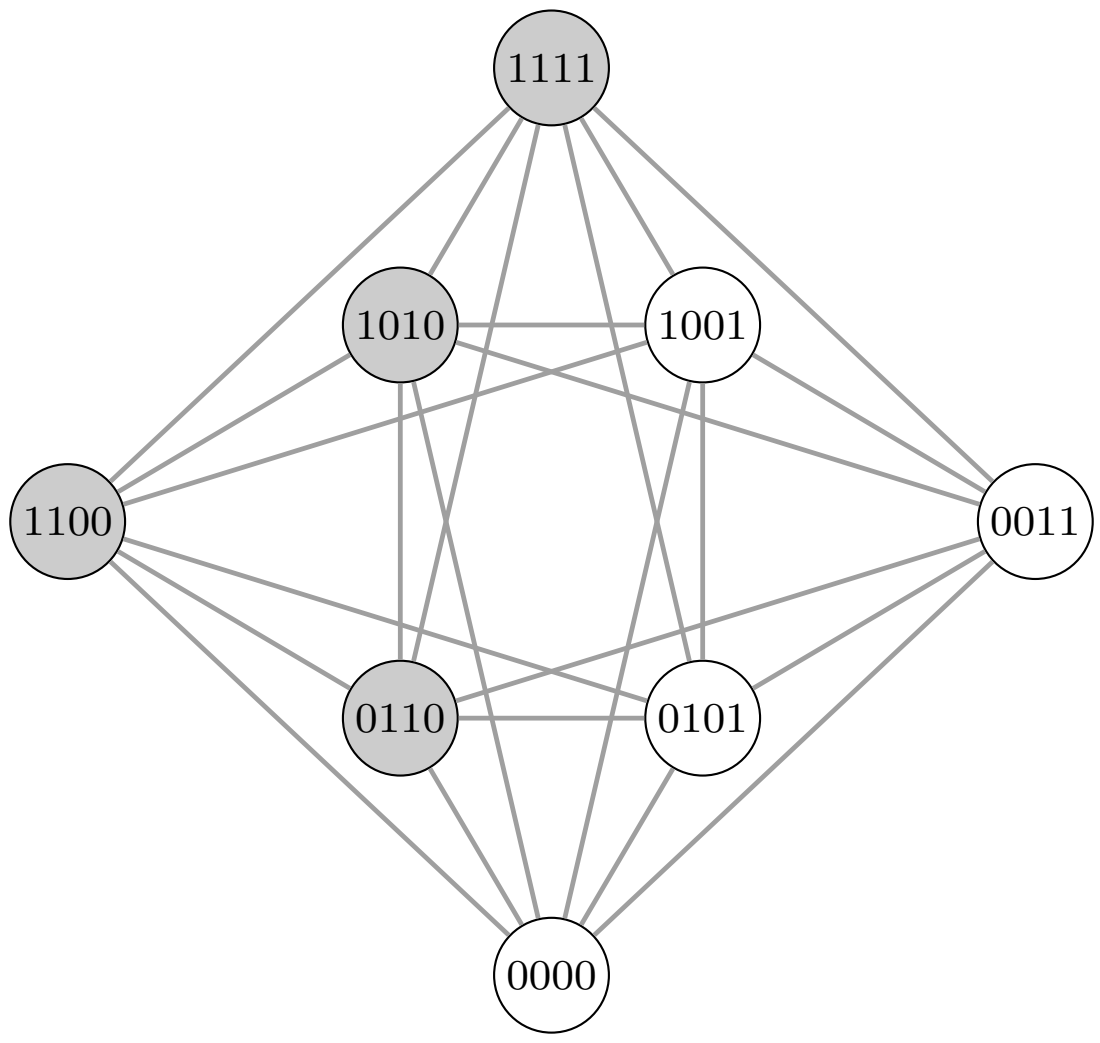
\begin{figure}[H]

\centering
\begin{tikzpicture}[
    scale=1.0,
    v/.style={circle,draw,inner sep=1.5pt,minimum size=18pt,font=\scriptsize,fill=white},
    code/.style={circle,draw,inner sep=1.5pt,minimum size=18pt,font=\scriptsize,fill=black!20},
    edge/.style={gray!75,line width=0.9pt}
]

\node[v]    (u0) at (0,-3.0) {0000};
\node[code] (u4) at (0, 3.0) {1111};

\node[code] (v1) at (-3.2,  0.0) {1100};
\node[code] (v2) at (-1.0,  1.3) {1010};
\node[v]    (v3) at ( 1.0,  1.3) {1001};
\node[code] (v4) at (-1.0, -1.3) {0110};
\node[v]    (v5) at ( 1.0, -1.3) {0101};
\node[v]    (v6) at ( 3.2,  0.0) {0011};

\foreach \x in {v1,v2,v3,v4,v5,v6}{
  \draw[edge] (u0)--(\x);
  \draw[edge] (u4)--(\x);
}

\draw[edge] (v1)--(v2);
\draw[edge] (v1)--(v3);
\draw[edge] (v1)--(v4);
\draw[edge] (v1)--(v5);

\draw[edge] (v6)--(v2);
\draw[edge] (v6)--(v3);
\draw[edge] (v6)--(v4);
\draw[edge] (v6)--(v5);

\draw[edge] (v2)--(v3);
\draw[edge] (v2)--(v4);
\draw[edge] (v3)--(v5);
\draw[edge] (v4)--(v5);

\end{tikzpicture}
\caption{The halved cube $Q_4$. The shaded vertices form the support
of $f^+$ in the $d=1$ case of Proposition~\ref{prop:sharp}.}
\label{fig:Q4-code}
\end{figure}

We note that the Boolean function in
Proposition~\ref{prop:sharp} defines a completely regular
code in $Q_{4d}$, obtained from the corresponding code
in $Q_4$ by the $\times d$ construction
\cite[Section~4.4]{KrotovMogilnykhVasileva2018}.
The code in $Q_4$ is depicted in
Figure~\ref{fig:Q4-code}.

\section{Completely regular codes of fixed degree in a halved cube}
\label{sec:codes}
A code in a graph is a subset of its vertex set.

 If $C$ is a code in $H_n$, put
\[
 C^+:=C\cap A_n.
\]
Let $C$ be a code in a graph $\Gamma$. Define
its covering radius $\rho_\Gamma(C)$ and distance layers by
\[
 \rho_\Gamma(C)=\max_{x\in V(\Gamma)}d_\Gamma(x,C),\qquad
 C_i=\{x\in V(\Gamma):d_\Gamma(x,C)=i\},\quad 0\le i\le\rho_\Gamma(C).
\]
The code is \emph{completely regular} (abbreviated \emph{CR})
if, for each $i$, every vertex of $C_i$ has the same numbers
$b_i$ and $c_i$ of neighbours in $C_{i+1}$ and $C_{i-1}$,
respectively. The intersection array is defined as
\[
 \iota_\Gamma(C)=
 \Iarr{b_0,\ldots,b_{\rho_{\Gamma}(C)-1}}{c_1,\ldots,c_{\rho_{\Gamma}(C)}}.
\]

For a nonempty code $C$ in $H_n$, let
\[
 \delta_H(C)=\min_{c\in C}
 |\{z\in C:d_{H_n}(c,z)=1\}|.
\]

For a real-valued function $f$ on a distance-regular graph $\Gamma$,
define its spectral support and external distance by
\[
I_\Gamma(f)=\{i:E_i(\Gamma)f\ne0\},\qquad
s_\Gamma(f)=|I_\Gamma(f)\setminus\{0\}|.
\]

A code $C$ is said to have degree $d$ if its characteristic
function $\bbone_C$ has degree $d$.
We further abbreviate
$I_\Gamma(C)=I_\Gamma(\bbone_C)$ and
$s_\Gamma(C)=s_\Gamma(\bbone_C)$.

Subscripts $H$ and $Q$ mean that the quantities are taken in
$H_n$ and $Q_n$, respectively.

\begin{theorem}\label{thm:external-distance}\cite[Theorems 5.10 and Theorem 5.13]{D73}
For a code $C$ in a distance-regular graph
$\Gamma$, one has the following 

1. $
 \rho_\Gamma(C)\le s_\Gamma(C).$
 
2. If $C$ is completely regular, then equality holds.

3. $s_\Gamma(C)=1$ if and only if $C$ is a CRC with $\rho_\Gamma(C)=1$.
\end{theorem}

\begin{proposition}\label{prop:distance}
If $C$ is a code in $H_n$ with $\delta_H(C)\ge1$, then
\[
 \rho_Q(C^+)\le
 \left\lceil \frac{\rho_H(C)}{2}\right\rceil.
\]
\end{proposition}
\begin{proof}
We show that  $d_{Q}(x,C^+)
 =\left\lceil \frac{d_{H}(x,C)}{2}\right\rceil$ and the bound for covering radius will follow. Fix $x\in A_n$, let $t=d_{H_n}(x,C)$, and choose a closest
vertex $c\in C$. Since $d_Q(x,y)=d_H(x,y)/2$ for every $y\in A_n$, we have
$d_Q(x,C^+)\ge\lceil t/2\rceil$.
If $t$ is even, then $c\in C^+$ and equality follows.
If $t$ is odd, then by $\delta_H(C)\ge1$,
there exists a neighbour $c'\in C$ of $c$ in $H_n$. It is of even weight and minimality of $t$ forces
$d_{H}(x,c')=t+1$, excluding $t-1$. Then
$d_{Q}(x,c')=(t+1)/2=\lceil t/2\rceil$, as required.
\end{proof}

\begin{theorem}\label{prop:main}
Let $C$ be a code in $H_n$, with
$\delta_H(C)\ge1$, and $$ i+j\ne n\qquad\text{for all }i,j\in I_H(C).$$
 Then
\[
 {\ C^+\text{ is CR in }Q_n
 \mbox{ if and only if }
 C\text{ is CR of covering radius }1\text{ in }H_n.\ }
\]
If either condition holds, then $\rho_Q(C^+)=1$.
\end{theorem}
\begin{proof}
By~\eqref{e} and the spectral condition,
$s_Q(C^+)=s_H(C)$.
If $C^+$ is completely regular, then
Theorem~\ref{thm:external-distance} and
Proposition~\ref{prop:distance} give
\[
 \rho_H(C)\le s_H(C)=s_Q(C^+)=\rho_Q(C^+)
 \le\left\lceil \frac{\rho_H(C)}{2}\right\rceil.
\]
Hence $\rho_H(C)=s_H(C)=1$.
Conversely, if $C$ is CR of radius $1$, then $s_Q(C^+)=s_H(C)=1$,
so $C^+$ is CR of covering radius $1$.
\end{proof}
Let $R_d$ be the maximum number of relevant variables that a degree-$d$ Boolean function on $H_n$ can have.
Recall that $R_d\leq \min(d\cdot 2^{d-1},4.394\cdot 2^{d})$ (see \cite{NS94,W22}).

\begin{corollary}\label{cor:large-n}
Let $d\ge1$ and suppose that
\[
 n\ge\max\{R_d,4d\}+1.
\]
Then every CRC $D$ of degree  $d$ in $Q_n$ has covering radius $1$ and the restriction $^+$ is a bijection between CRCs of covering radius $1$ of degree $d$ on $H_n$ and $Q_n$.
\end{corollary}
\begin{proof}
Let $D$ be completely regular of degree $d$ in $Q_n$,
and let $g=\bbone_D$. By Theorem~\ref{cor:boolean-lifting},
$g$ has a unique Boolean extension $f$ of degree  $d$
to $H_n$. Set $C$ as follows: $\bbone_C=f$, $C^+=D$.
Since $n\geq R_d+1$, $f$ has at least one irrelevant variable.
Therefore, we have $\delta_H(C)\ge1$. 
Since $I_H(C)\subseteq\{0,\ldots,d\}$
and $n>2d$, Proposition~\ref{prop:main} applies.
The converse follows in the same way as at the end of the proof of Theorem  \ref{prop:main}.

\end{proof}

{\bf Remark 1.}
In contrast to halved cubes, the classical Hamming, Johnson, and Grassmann graphs admit completely regular codes of degree $d$ and covering radius $d$ for every admissible $d$. A simple construction is to fix an $(n-d)$-subset $S$ and take all vertices whose supports are contained in $S$. This works for both Hamming and Johnson graphs.
In the latter case, vertices are binary vectors of
weight $k$, and we assume $d\le k\le n/2$. In Grassmann graphs, the analogous construction is obtained by replacing subsets with subspaces.  By Corollary~\ref{cor:large-n}, however, for every fixed $d$, all CRCs of degree $d$ in sufficiently large halved cubes have covering radius $1$.
\section{Functions of degree at most two on $Q_n$}
\label{sec:classification}

By Theorem~\ref{cor:boolean-lifting}, for $n\ge5$,
every Boolean function of degree $1$ on
$Q_n$ extends uniquely to a Boolean function of degree
$1$ on $H_n$ and is thus equivalent to $x_1^+$.
For $n=4$, there are exactly two equivalence classes
under automorphisms induced by $H_4$ and complementation:
$x_1^+$ and the restriction of the function from
Proposition~\ref{prop:sharp} with $d=1$. We leave the verification of this classification to the reader.

For $n\ge9$, Theorem~\ref{cor:boolean-lifting} reduces the
classification of Boolean functions of degree $2$ on $Q_n$
to that on $H_n$. Up to automorphisms and complementation,
the five representatives are listed in Table~\ref{tab:boolean-reps}.
Here $m$ is the number of relevant coordinates; the others are
irrelevant. By Corollary~\ref{cor:large-n}, every CRC of degree $2$ in $Q_n$
for $n\ge9$ is equivalent 
to one of the codes $C^+$ marked CR in the table. In particular, our results settle the two open cases
for the eigenvalue $\theta_2(12)$
in \cite[Section~6.5]{KrotovMogilnykhVasileva2018}.

Let
\[
 S_8=\{0000,1000,1100,1110,1111,0111,0011,0001\}.
\]
Here $S_8$ induces an $8$-cycle in $H_4$, and
$\bbone_{S_8}(x_1,x_2,x_3,x_4)$ means the characteristic function of that
set. In row 4 the function has support
$\{010,011,101,111\}\subseteq Z_2^3$.

\begin{table}[H]
\centering
\footnotesize
\renewcommand{\arraystretch}{1.32}
\setlength{\tabcolsep}{1.2pt}
\begin{tabularx}{0.73\linewidth}{@{}
  >{\raggedright\arraybackslash}X
  c @{\hspace{8pt}}c c
  >{\centering\arraybackslash}p{1.5cm}
  c
  >{\centering\arraybackslash}p{2.65cm}
  @{}}
\toprule
 & & & \multicolumn{2}{c}{Code $C$ in $H_n$}
       & \multicolumn{2}{c}{Code $C^+$ in $Q_n$} \\
\cmidrule(lr){4-5}\cmidrule(l){6-7}
Representative $f$ & $m$ & $I_H(f)\setminus \{0\}$
 & CR & $\iota_H(C)$ & CR & $\iota_Q(C^+)$ \\
\midrule

$f=\bbone_{\{00,11\}}(x_1,x_2)$ & 2 & $\{2\}$
 & $\yes$ & $\Iarr{2}{2}$
 & $\yes$ & $\Iarr{2(n-2)}{2(n-2)}$ \\
$f=x_1x_2$ & 2 & $\{1,2\}$
 & $\yes$ & $\Iarr{2,1}{1,2}$
 & $\no$ & -- \\
$f=\bbone_{\{000,111\}}(x_1,x_2,x_3)$ & 3 & $\{2\}$
 & $\yes$ & $\Iarr{3}{1}$
 & $\yes$ & $\Iarr{3(n-2)}{n-2}$ \\
$f=(1-x_1)x_2+x_1x_3$ & 3 & $\{1,2\}$
 & $\no$ & --
 & $\no$ & -- \\
$f=\bbone_{S_8}(x_1,x_2,x_3,x_4)$ & 4 & $\{2\}$
 & $\yes$ & $\Iarr{2}{2}$
 & $\yes$ & $\Iarr{2(n-2)}{2(n-2)}$ \\
\bottomrule
\end{tabularx}
\caption{Representatives of Boolean functions of degree 
$2$ on $H_n$, $n\ge9$, up to hypercube automorphisms and
complementation $f\mapsto1-f$. Here $C$ is the support of $f$;
the code parameters refer to the displayed representatives.}
\label{tab:boolean-reps}
\end{table}

We also note that in the same manner degree $3$ Boolean functions on $Q_n$ from degree $3$ Boolean functions on $H_n$ (classified in \cite[Section 2.1]{KV24}) could be described.

{\bf Acknowledgments.} A key observation used in the proof of Theorem \ref{cor:boolean-lifting} emerged during a discussion of a related problem with the AI system Aristotle. The idea was subsequently adapted by the authors to the present setting.

\end{document}